\documentclass[10pt,a4paper]{amsart}

\begin{document}

\newtheorem{lemma}{Lemma}[section]
\newtheorem{prop}[lemma]{Proposition}
\newtheorem{cor}[lemma]{Corollary}
\newtheorem{thm}[lemma]{Theorem}

\theoremstyle{definition}
\newtheorem{rem}[lemma]{Remark}
\newtheorem{rems}[lemma]{Remarks}
\newtheorem{defi}[lemma]{Definition}
\newtheorem{ex}[lemma]{Example}

\newcommand{\slth}{SL_3({\bf H})}
\newcommand{\slto}{E_{6(-26)}}
\newcommand{\sltc}{SL_3({\bf C})}
\newcommand{\suto}{F_{4(-52)}}
\newcommand{\sltr}{SL_3({\bf R})}
\newcommand{\sltk}{SL_3({\bf K})}
\newcommand{\sutk}{SU_3({\bf K})}
\newcommand{\sutc}{SU_3({\bf C})}
\newcommand{\suth}{SU_3({\bf H})}
\newcommand{\sutr}{SO_3({\bf R})}
\title[Finite products of hyperbolic spaces with $SL_3$'s]{Property (RD) for cocompact lattices in a finite product of rank one Lie groups with some rank two Lie groups}
\author{Indira Chatterji}
\address{ETHZ, D-MATH, R\"amistra\ss e 101, 8092 Z\"urich}
\email{indira@math.ethz.ch}
\date{April 2001}
\begin{abstract}We apply V. Lafforgue's techniques to establish property (RD) for cocompact lattices in a finite product of rank one Lie groups with Lie groups whose restricted root system is of type $A_2$.\end{abstract}
\maketitle
We recall that a {\it length function} on a discrete group $\Gamma$ is a function $\ell:\Gamma\to{\bf R}_+$ such that the neutral element is mapped to zero, and such that $\ell(\gamma)=\ell(\gamma^{-1})$, $\ell(\gamma\mu)\leq\ell(\gamma)+\ell(\mu)$ for any $\gamma,\mu\in\Gamma$. A discrete group $\Gamma$ is said to have {\it property (RD) with respect to $\ell$} if there exists a polynomial $P$ such that for any $r\in{\bf R}_+$ and $f\in{\bf C}\Gamma$ supported on elements of length shorter than $r$ the following inequality holds:
$$\|f\|_*\leq P(r)\|f\|_2$$
where $\|f\|_*$ denotes the operator norm of $f$ acting by left convolution on $\ell^2(\Gamma)$, and $\|f\|_2$ the usual $\ell^2$ norm. Property (RD) has been first established for free groups by Haagerup in \cite{Haagerup}, but introduced and studied by P. Jolissaint in \cite{Jolissaint}, who established it for classical hyperbolic groups. The extension to Gromov hyperbolic groups is due to P. de la Harpe in \cite{Harpe}. Providing the first examples of higher rank groups, J. Ramagge, G. Robertson and T. Steger in \cite{RRS} proved that property (RD) holds for discrete groups acting freely on the vertices of an $\tilde{A}_1\times\tilde{A}_1$ or $\tilde{A}_2$ building and recently V. Lafforgue did it for cocompact lattices in $\sltr$ and $\sltc$ in \cite{Lafforgue}. A conjecture, due to A. Valette, see \cite{Ferry} claims that property (RD) holds for any discrete group acting isometrically, properly and cocompactly either on a Riemannian symmetric space or on an affine building. Property (RD) is important in the context of Baum-Connes conjecture, precisely, V. Lafforgue in \cite{Lafforgue2} proved that for ``good'' groups having property (RD), the Baum-Connes conjecture without coefficients holds.

In this article we establish property (RD) for discrete cocompact subgroups $\Gamma$ of a finite product of type
$${\rm Iso}({\mathcal X}_1)\times\dots\times{\rm Iso}({\mathcal X_n})$$
where the ${\mathcal X}_i$'s are either Gromov hyperbolic spaces, buildings associated to $SL_3(F)$ (for $F$ a non-archimedian locally compact field), or symmetric spaces associated to Lie groups whose restricted root system is of type $A_2$ (those are known to be locally isomorphic to $\sltr$, $\sltc$, $\slth$ and $\slto$, see \cite{Helgason}). In particular, this immediately implies the following
\begin{thm}\label{principal}
Any cocompact lattice in
$$G=G_1\times\dots\times G_n$$
has property (RD), where the $G_i$'s are either rank one Lie groups, $\sltr$, $\sltc$, $\slth$ or $\slto$.
\end{thm} 
To establish this result, we will first answer (positively) a question posed by V. Lafforgue in \cite{Lafforgue}, which was to know whether his Lemmas 3.5 and 3.7 are still true for the groups $\slth$ and $\slto$, whose associated symmetric spaces have also flats of type $A_2$. Observing that these lemmas are in fact ``three points conditions'' it will be enough to prove that if $X$ denotes $\slth/\suth$ or $\slto/\suto$ then for any three points in $X$ there exists a totally geodesic embedding of $\sltc/\sutc$ containing those three points. Secondly we will explain how to use the techniques used in \cite{RRS} and \cite{Lafforgue} for the above described products. This provides many new examples of groups with property (RD), such as irreducible cocompact lattices in $SL_2({\bf R})\times SL_2({\bf R})$.

Finally we remark that combining our result with V. Lafforgue's crucial theorem in \cite{Lafforgue2} yields the following
\begin{cor}
The Baum-Connes conjecture without coefficients (see \cite{Valette}) holds for any cocompact lattice in
$$G=G_1\times\dots\times G_n$$
where the $G_i$'s are either rank one Lie groups, $\sltr$, $\sltc$, $\slth$ or $\slto$.
\end{cor}
It is a pleasure to thank my thesis advisor A. Valette for suggesting me this problem and for all his generous help. I must also thank V. Lafforgue for graciously explaining the methods of \cite{RRS} and \cite{Lafforgue} to me. I am also extremely indebted to D. Allcock and G. Prasad for important discussions on the exceptional group $\slto$. 
%
\section{The case of $\slth$.}
We will write ${\bf H}$ for {\it Hamilton's Quaternion algebra}, which is a 4 dimensional vector space over ${\bf R}$, whose basis is given by the elements $1,i,j$ and $k$, satisfying
$$i^2=j^2=k^2=-1\,,\ ij=k\,,\ ki=j\,,\ jk=i.$$
This is an associative division algebra endowed with an involution $h\mapsto\overline{h}$ which is the identity over $1$, and minus the identity over $i,j$ and $k$. A norm on ${\bf H}$ can be given by $|h|=\sqrt{h\overline{h}}\in{\bf R}_+$. A quaternion will be called a {\it unit} if of norm one, {\it real} if lying in ${\rm span}\{1\}$ and {\it imaginary} if lying in ${\rm span}\{1\}^\bot$ (for the scalar product of ${\bf R}^4$ which turns the above described basis in an orthonormal basis). Note that an imaginary unit has square $-1$. The following lemma is obvious:
\begin{lemma}\label{C_dans_H} Any element $h\in{\bf H}$ is contained in a commutative subfield of ${\bf H}$.
\end{lemma}
\hfill\qedsymbol
\begin{defi}
Denote by $I$ the identity matrix in $M_3({\bf C})$, and
$$J=\left(\begin{array}{cc}0&-I\\
I&0\end{array}\right)\in M_6({\bf C}),$$
we write 
$$M_3({\bf H})=\{M\in M_6({\bf C}) \hbox{ such that } JMJ^{-1}=\overline{M}\}$$
$$\slth=\{M\in M_3({\bf H}) \hbox{ such that } {\rm det}(M)=1\}$$
$$\suth=\{M\in\slth \hbox{ such that }MM^*=M^*M=I\}.$$
\end{defi}
\begin{rem}\label{decomposition_de_H} 1) Given any imaginary unit $\mu\in{\bf H}$, we can write ${\bf H}={\bf C}\oplus{\bf C}\mu$, and thus decompose any $h\in{\bf H}$ as $h=h_1+h_2\mu$, with $h_1,h_2\in{\bf C}$. Similarly, any $3\times 3$ matrix $M$ with coefficients in ${\bf H}$ can be written $M=M_1+M_2\mu$, with $M_1,M_2\in M_3({\bf C})$. The *-homomorphism
$$M\mapsto\left(\begin{array}{cc}M_1&-\overline{M_2}\\
M_2&\overline{M_1}\end{array}\right)$$
gives then an isomorphism between the algebra of $3\times 3$ matrices with coefficients in ${\bf H}$ (and usual multiplication) and $M_3({\bf H})$ as just defined. It is a straightforward computation to see that elements of type
$${\mathcal T}_\lambda=\left(\begin{array}{ccc}\lambda&0&0\\
0&1&0\\
0&0&1\end{array}\right)$$
with $\lambda$ a unit in ${\bf H}$ actually belong to $\suth$.
\end{rem}
\begin{defi} Let ${\bf K}$ denote ${\bf R},{\bf C}$ or ${\bf H}$. We set
$$X_{\bf K}=SL_3({\bf K})/SU_3({\bf K}).$$
It will be convenient for us to choose
$$\{M\in SL_3({\bf K})\hbox{ such that } M^*=M, M\hbox{ positive}\}$$
as a model for $X_{\bf K}$. On $X_{\bf K}$ we consider the action of $\sltk$ given by:
\begin{eqnarray*}\sltk\times X_{\bf K}&\to & X_{\bf K}\\
(\,g\ ,\,z\,) &\mapsto & g(z)=(gz^2g^*)^{1/2}.
\end{eqnarray*}
This action is transitive since for $M\in X_{\bf K}$, setting $g=M$ we get that $g(I)=M$. Moreover the stabilizer of $I$ is clearly $\sutk$. We equip $X_{\bf K}$ with the distance
$$d_{\bf K}(x,y)=\log\|x^{-1}y\|+\log\|y^{-1}x\|.$$
where $\|\ \|$ denotes the operator norm on $\sltk$ acting on ${\bf K}^3$. Notice that for $a\in X_{\bf K}$ a diagonal matrix (thus real), if we assume that
$$a=\left(\begin{array}{ccc}e^{\alpha_1}&0&0\\
0&e^{\alpha_2}&0\\
0&0&e^{\alpha_3}\end{array}\right)\hbox{ with }\alpha_1\geq\alpha_2\geq\alpha_3\hbox{ and }\alpha_1+\alpha_2+\alpha_3=0$$
we get that $d_{\bf K}(a,I)=\alpha_1-\alpha_3$. Denote by $A$ the set of all diagonal matrices in $X_{\bf K}$. To stick with the classical terminology, we will call {\it flats} the sets of type $g(A)$ for a $g\in\sltk$. 
\end{defi}
\begin{rem}\label{generateurs_de_slth}  
The action of $\sltk$ on $X_{\bf K}$ is by isometries with respect to the above given distance (in other words, the distance is $\sltk$-invariant): For $x,y\in X_{\bf K}$, we have to show that the operator norm of $x^{-1}y$ is equal to that of $z=g(y)^{-1}g(x)=(gx^2g^*)^{-1/2}(gy^2g^*)^{1/2}$. But for any $k_1,k_2\in\sutk$, $\|z\|=\|k_1zk_2\|$, so setting $k_1=(gx)^{-1}(gx^2g^*)^{1/2}$ and $k_2=(gy^2g^*)^{-1/2}gy$, we get that $k_1zk_2=x^{-1}y$. Notice that for $g\in\sltk$, the standard action $z\mapsto gzg^*$ is not isometric. Indeed, the operator norm is in general not invariant by conjugation by $g$ unless $g=k\in\sutk$, and in that case $k(z)=(kz^2k^*)^{1/2}=(kzk^*kzk^*)^{1/2}=kzk^*$.
\end{rem}
\begin{lemma}\label{plongeCdansH} Let $\varphi:{\bf C}\to{\bf H}$ be an isometric injective ring homomorphism. It induces a group homomorphism
$$\varphi^*:\sltc\to\slth$$
which induces a totally geodesic isometric embedding
$$\overline{\varphi}:X_{\bf C}\to X_{\bf H}.$$
\end{lemma}
\begin{proof}
That $\varphi^*$ is a homomorphism is clear in view of Remark \ref{decomposition_de_H}, so let us now show that $\overline{\varphi}$ is an isometry. Since $\sltc$ acts transitively and by isometries on $X_{\bf C}$, it will be enough to show that 
$$d_{\bf H}(\overline{\varphi}(x),I)=d_{\bf C}(x,I)$$ 
for any $x\in X_{\bf C}$. Noticing that $\varphi^*(\sutc)\subset\suth$ we deduce that for any $g\in\sltc$, $\overline{\varphi}\left(g(x)\right)=\varphi^*(g)\left(\overline{\varphi}(x)\right)$. Now, for $x\in X_{\bf C}$, there exists $k\in\sutc$ such that $k(x)=a$ is diagonal with positive coefficients. Since $\varphi$ is a ring homomorphism, it is the identity over the real numbers and so:
\begin{eqnarray*}d_{\bf C}(x,I) &=& d_{\bf C}(k(x),k(I))=d_{\bf C}(a,I)=d_{\bf H}(a,I)\\
&=&d_{\bf H}(\varphi^*(k^*)(a),\varphi^*(k^*)(I))=d_{\bf H}(\overline{\varphi}(k^*(a)),I)\\
&=&d_{\bf H}(\overline{\varphi}(x),I).
\end{eqnarray*}
Let us now prove that the embedding is totally geodesic. Take $x,y$ in the image of $\overline{\varphi}$, so that $x=\overline{\varphi}(x')$ and $y=\overline{\varphi}(y')$ for $x',y'\in X_{\bf C}$. There is a $g\in\sltc$ so that $x'=g(I)$ and $y'=g(a)$ for $a\in A$. Then $x= \overline{\varphi}(g(I))=\varphi^*(g)(I)$ and $y=\overline{\varphi}(g(a))=\varphi^*(g)(a)$. Since $A$ is obviously in the image of $\overline{\varphi}$, so will $\varphi^*(g)(A)$ be (that is, the whole flat containing $x$ and $y$ is in the image of the embedding). Since a geodesic $\gamma$ between $x$ and $y$ lies in any flat containing them, we conclude that $\gamma\subset\varphi^*(g)(A)\subset\overline{\varphi}(X_{\bf C})$. 
\end{proof}
\begin{rem} Any element of $X_{\bf H}$ can be diagonalized using elements of type ${\mathcal T}_h$ (described in Remark \ref{decomposition_de_H}) and of $\sutr$. Indeed, take $z=(z_{ij})\in X_{\bf H}$, then $z'={\mathcal T}_h(z)$ with $h=\overline{z_{12}}/|z_{12}|$ verifies that $z'_{ij}\in{\bf R}$ for $i,j=1,2$ (in case $z_{12}=0$, or if it lies in ${\bf R}$ we can just skip this step). Take $k\in SO_2({\bf R})$ diagonalizing the $2\times 2$ matrix $(z'_{ij})$ $i,j=1,2$. Then
\begin{eqnarray*}z''& = & \left(\begin{array}{ccc}
0&0&1\\
1&0&0\\
0&1&0\end{array}\right)
\left(\begin{array}{cc}
k&0\\
0&1\end{array}\right)z'\left(\begin{array}{cc}
k^*&0\\
0&1\end{array}\right)
\left(\begin{array}{ccc}
0&1&0\\
0&0&1\\
1&0&0\end{array}\right)\\
&=&\left(\begin{array}{ccc}
z''_{11}&z''_{12}&z''_{13}\\
\overline{z''_{12}}&z''_{22}&0\\
\overline{z''_{13}}&0&z''_{33}\end{array}\right)\end{eqnarray*}
with $z''_{11},z''_{22},z''_{33}\in{\bf R}$. Repeating this argument twice we get the diagonalization.

\smallskip

The above given argument shows that the subgroup of $\slth$ generated by the ${\mathcal T}_h$'s and $\sltr$ acts transitively on $X_{\bf H}$. Since the ${\mathcal T}_h$'s and $\sutr$ do stabilize $I$, we can conclude that $\slth$ is generated by elements of type ${\mathcal T}_h$ and $\sltr$.
\end{rem}
\begin{prop}\label{principale_h}For any three points in $X_{\bf H}$ there is a totally geodesic embedding of $X_{\bf C}$ containing those three points.\end{prop}
\begin{proof}Let $x,y,z\in X_{\bf H}$. Up to multiplication by an element of $\slth$ (which is an isometry), we can assume that $x=I$ and $y$ is a diagonal matrix (see the previous Remark). We write
$$z=\left(\begin{array}{ccc}z_{11} & z_{12} & z_{13}\\
\overline{z_{12}} & z_{22} & z_{23}\\
\overline{z_{13}} & \overline{z_{23}} & z_{33}\end{array}\right)$$
with $z_{11},z_{22},z_{33}\in{\bf R}$ and $z_{12},z_{13}$ and $z_{23}$ in ${\bf H}$. Consider the element
$$k=\left(\begin{array}{ccc}\overline{z_{12}}/|z_{12}|&0&0\\
0&1&0\\
0&0&{z_{23}}/|z_{23}|\end{array}\right)\in\suth,$$
then $k(z)$ only has real elements except for $k(z)_{13}=\overline{k(z)_{31}}$, and $k(a)=a$ for any diagonal element $a$. Applying Lemma \ref{C_dans_H} to $k(z)_{13}$, we get an embedding of ${\bf C}$ in ${\bf H}$ which contains $k(z)_{13}$, and thus a totally geodesical embedding $\overline{\varphi}:X_{\bf C}\to X_{\bf H}$ whose image contains $x,y$ and $z$.
\end{proof}
\section{The case of $\slto$.}
We will write ${\bf O}$ for the 8 dimensional algebra over ${\bf R}$, whose basis is given by the elements $e_0,e_1,\dots,e_7$, satisfying
\begin{eqnarray*}& e_ie_0=e_0e_i=e_i, e_i^2=-e_0 & \hbox{ for all }i=1,\dots,7\\
&e_ie_j=-e_je_i &\hbox{ if }i\not= j \hbox{ and } i,j\not= 0\\
&e_2e_6=e_3e_4=e_5e_7=e_1
\end{eqnarray*}
and all those one can deduce by cyclically permuting the indices from 1 to 7. This is a non-associative division algebra, endowed with an involution $x\mapsto\overline{x}$ which is the identity over $e_0$, and minus the identity over $e_i$, for all $i=1,\dots,7$. A norm on ${\bf O}$ can be given by $|x|=\sqrt{x\overline{x}}\in{\bf R}_+$. The division algebra ${\bf O}$ is called the {\it Cayley Octonions}. An octonion will be called a {\it unit} if of norm one, {\it real} if lying in ${\rm span}\{e_0\}$ and {\it imaginary} if lying in ${\rm span}\{e_0\}^\bot$ (for the scalar product of ${\bf R}^8$ which turns the above described basis in an orthonormal basis). The following
\begin{thm}[Artin, see \cite{Sch} page 29]\label{H_dans_O} Any two elements of ${\bf O}$ are contained in an associative sub-division algebra of ${\bf O}$.
\end{thm}
\noindent
will allow us to apply the arguments used in the previous section to $\slto$.

We will now give a definition of $\slto$, of a maximal compact subgroup $\suto$, and of a model for $X_{\bf O}=\slto/\suto$. This part is based on the work of H. Freudenthal, see \cite{Freud1} and \cite{Freud2}.
\begin{defi}Denote by $M_3({\bf O})$ the set of $3\times 3$ matrices with coefficients on ${\bf O}$. For $M=(m_{ij})\in M_3({\bf O})$ write $M^*=\overline{M}^t=(\overline{m}_{ji})$. The {\it exceptional Jordan algebra} is given by
$${\mathcal J}=\{M\in M_3({\bf O})\hbox{ such that }M=M^*\}$$
which is stable under the {\it Jordan multiplication} given by
$$M\star N={\frac {1}{2}}(MN+NM).$$
H. Freudenthal defined an application
\begin{eqnarray*}{\det}:{\mathcal J} & \to & {\bf R}\\
\left(\begin{array}{ccc}\xi_1&x_3&\overline{x}_2\\
\overline{x}_3&\xi_2&x_1\\
x_2&\overline{x}_1&\xi_3\end{array}\right)& \mapsto & \xi_1\xi_2\xi_3-(\sum_{i=1}^3\xi_i|x_i|^2)+2\Re (x_1x_2x_3)
\end{eqnarray*}
(where $\Re$ denotes the real part, and is well defined even without parenthesis), showed that 
$$\slto=\{g\in GL({\mathcal J})\hbox{ such that }\det\circ g=\det\},$$
and that
$$\suto=\{g\in\slto\hbox{ such that }g(I)=I\}$$
is the automorphism group of ${\mathcal J}$ and is a maximal compact subgroup in $\slto$.
\end{defi}
\begin{rems}\label{elements_de_slto} We will now explicitly show some elements of $\slto$ (see \cite{Freud2} for the proofs of their belonging to $\slto$):

1) Any element $x$ of $\sltr$ gives a map $x:{\mathcal J}\to{\mathcal J}$ by $M\mapsto x(M)=xMx^t$ which preserves the determinant.

2) Let $a$ be a unit in ${\bf O}$, define 
$$\varphi_a=\left(\begin{array}{ccc}1&0&0\\
0&a&0\\
0&0&\overline{a}\end{array}\right).$$
The map $\psi_a:{\mathcal J}\to{\mathcal J}$, defined as $\psi_a(M)=\varphi_{\overline{a}}M\varphi_a$ (and no parenthesis are needed) is in $\slto$, and even in $\suto$.

3) H. Freudenthal proves (see \cite{Freud1}, page 40) that elements in ${\mathcal J}$ are diagonalizable by elements in $\suto$, that the elements on the diagonal are uniquely determined up to permutation.
\end{rems}
\begin{defi}
We say that an element $M\in{\mathcal J}$ is {\it positive} if after diagonalization it only has positive elements. We define
$$X_{\bf O}=\{M\in {\mathcal J}\hbox{ such that }\det (M)=1 \hbox{ and positive}\}.$$
We let $\slto$ act on $X_{\bf O}$ as follows:
\begin{eqnarray*}\slto\times X_{\bf O}&\to&X_{\bf O}\\
(g,M)&\mapsto&g\circ M=\sqrt{g(M^2)}
\end{eqnarray*}
where (for $M$ positive) $\sqrt{M}$ is an element $Y$ in $X_{\bf O}$ satisfying $Y\star Y=M$, and well defined by requiring positivity. Notice that for $k\in\suto$, we have that $k\circ M=k(M)$.
\end{defi}
\begin{prop}The group $\slto$ acts transitively on $X_{\bf O}$, and the stabilizer of $I$ is $\suto$.\end{prop}
\begin{proof}Since $k\circ I=k(I)=I$, the stabilizer of $I$ is $\suto$, and the action is transitive because given any $M\in X_{\bf O}$, there exists a $k\in\suto$ such that $k\circ M=D$ is diagonal. Since $\det (M)=\det\left(k(M)\right)=1$, the element $D^{-1}$ defines an element in $\slto$ (coming from $SL_3({\bf  R})$), and $D^{-1}\circ D=I$.
\end{proof}
\begin{rem}\label{deux_elements_dans_un_plat} For any two elements $M,N$ of $X_{\bf O}$, there exists a $g\in\slto$ such that $M=g\circ I$ and $N=g\circ D$ where $D$ is a diagonal matrix in $X_{\bf O}$. Indeed, in the proof of the preceding proposition we saw that there exists an $h\in\slto$ with $I=h\circ M$, and that any matrix (and thus in particular $h\circ N$) is diagonalizable with elements of $\suto$. We then choose $k\in\suto$ such that $k(h\circ N)=D$ is diagonal, and set $g=(kh)^{-1}$.
\end{rem}
\begin{defi}For a diagonal element $D$ in $X_{\bf O}$, we set 
$$\|D\|=\max\{d_i|i=1,2,3\}$$
where the $d_i$'s denote the elements in the diagonal of $D$. We equip $X_{\bf O}$ with the distance 
$$d(M,N)=\log\|D\|+\log\|D^{-1}\|.$$ 
for $D$ as in the previous remark. It is well defined because of point 3) of Remark \ref{elements_de_slto}.
\end{defi}
Any embedding of ${\bf H}$ in ${\bf O}$ gives a decomposition ${\bf O}={\bf H}\oplus{\bf H}\ell$ (where $\ell$ is any unit in the orthogonal complement of the embedded copy of ${\bf H}$ in ${\bf O}$), and we will now describe some other elements in $\suto$ that will later be useful to embed $\slth$ in $\slto$ starting from an embedding of ${\bf H}$ in ${\bf O}$. Given a unit $h\in{\bf H}$ we define two maps
\begin{eqnarray*}t_h,u_h:{\bf O}&\to & {\bf O}\\
x=a+b\ell&\mapsto &t_h(x)=ha+b\ell\\
& & u_h(x)=a+(b\overline{h})\ell
\end{eqnarray*}
and an element (that we will later prove to belong to $\suto$)
\begin{eqnarray*}T_h:{\mathcal J}&\to & {\mathcal J}\\
\left(\begin{array}{ccc}
\xi_1&x_3&\overline{x}_2\\
\overline{x}_3&\xi_2&x_1\\
x_2&\overline{x}_1&\xi_3
\end{array}\right)
& \mapsto & 
\left(\begin{array}{ccc}
\xi_1&t_h(x_3)&t_h(\overline{x}_2)\\
\overline{t_h(x_3)}&\xi_2&\overline{u_h(\overline{x}_1)}\\
\overline{t_h(\overline{x}_2)}&u_h(\overline{x}_1)&\xi_3
\end{array}\right)
\end{eqnarray*}
Using that for $x=a+b\ell$ and $y=c+d\ell$, the product $xy$ reads $(ac-\overline{d}b)+(b\overline{c}+da)\ell$ (see \cite{Springer}), it is a direct computation to check that $\det\circ T_h=\det$.

Knowing that $T_h\in\slto$, it is pretty clear that $T_h(I)=I$, and thus $T_h\in\suto$.
We will now, for any given embedding of ${\bf H}$ in ${\bf O}$, define an explicit embedding of $\slth$ in $\slto$.
\begin{lemma}\label{plongeHdansO}
Let $S$ denote the subgroup of $\slto$ generated by elements of the type $x\in\sltr$ and $T_h$ for $h$ a unit in ${\bf H}$ as just described. Then $S$ is isomorphic to $\slth$.
\end{lemma}
\begin{proof}An embedding of ${\bf H}$ in ${\bf O}$ defines an embedding of $X_{\bf H}$ in ${\mathcal J}$ by taking the images in ${\bf O}$ of the coefficients of an element in $X_{\bf H}$. It will now be enough to see that the map
\begin{eqnarray*}\rho:S &\to &\slth\\
z &\mapsto & z|_{X_{\bf H}}
\end{eqnarray*}
is well defined and a group isomorphism.

Remember that for $h$ a unit in ${\bf H}$ we denoted by ${\mathcal T}_h\in\slth$ a diagonal matrix with $h$ in the first place and $1$'s in the second and third place. Now, notice (from the definition of ${\mathcal T}_h$) that $\rho(T_h)={\mathcal T}_h$ for any $h$ a unit in ${\bf H}$ and that trivially, for $x\in\sltr$, $\rho(x)\in\slth$. This shows that the image of $\rho$ is contained in $\slth$. The map $\rho$ being a restriction, it is a group homomorphism. 

Our map $\rho$ is clearly injective, and because of Remark \ref{generateurs_de_slth} we know that $\rho$ maps the generators of $S$ onto those of $\slth$.
\end{proof}
\begin{prop}\label{principale_o} For any three points in $X_{\bf O}$ there is a totally geodesic embedding of $X_{\bf C}$ containing those three points.\end{prop}
\begin{proof}Let $x,y,z\in X_{\bf O}$. Up to multiplication by an element of $\slto$ (which is an isometry), we can assume that $x=I$ and $y$ is a diagonal matrix (use Remark \ref{deux_elements_dans_un_plat}). We write
$$z=\left(\begin{array}{ccc}z_{11} & z_{12} & z_{13}\\
\overline{z_{12}} & z_{22} & z_{23}\\
\overline{z_{13}} & \overline{z_{23}} & z_{33}\end{array}\right)$$
with $z_{11},z_{22},z_{33}\in{\bf R}$ and $z_{12},z_{13}$ and $z_{23}$ in ${\bf O}$. Consider the element $k=\phi_a$ where $a=\overline{z_{12}}/|z_{12}|$ is a unit in ${\bf O}$
then $k(z)$ only has real elements except for $k(z)_{13}=\overline{k(z)_{31}}$, and $k(d)=d$ for any diagonal element in $X_{\bf O}$. Applying Theorem \ref{H_dans_O} to $k(z)_{13}$ and $k(z)_{23}$, we get an embedding of ${\bf H}$ in ${\bf O}$ which contains $k(z)_{13}$ and $k(z)_{23}$, and thus an embedding $\overline{\varphi}:X_{\bf H}\to X_{\bf O}$ whose image contains $x,y$ and $z$, and which is isometric by definition of the distance on $X_{\bf O}$.

To see that the embedding is totally geodesic is basically the same argument as in the proof of Proposition \ref{plongeCdansH}: Take $x,y$ in the image of $\overline{\varphi}$, so that $x=\overline{\varphi}(x')$ and $y=\overline{\varphi}(y')$ for $x',y'\in X_{\bf C}$. There is a $g\in\sltc$ so that $x'=g(I)$ and $y'=g(a)$ for $a\in A$. Now, the embedding $\overline{\varphi}$ comes from an embedding $\varphi*:\slth\to\slto$ described in Lemma \ref{plongeHdansO}, so that $x=\overline{\varphi}(g(I))=\varphi^*(g)(I)$ and $y=\overline{\varphi}(g(a))=\varphi^*(g)(a)$. Since $A$ is obviously in the image of $\overline{\varphi}$, so will $\varphi^*(g)(A)$ be (that is, the whole flat containing $x$ and $y$ is in the image of the embedding). Since a geodesic $\gamma$ between $x$ and $y$ lies in any flat containing them, we conclude that $\gamma\subset\varphi^*(g)(A)\subset\overline{\varphi}(X_{\bf C})$. 

We are now reduced to find an embedding of $X_{\bf C}$ in $X_{\bf H}$ containing $x,y$ and $z$, but we get this one because of Proposition \ref{principale_h}, and combining it with $\overline{\varphi}$ we get the sought embedding.
\end{proof}
\section{How we answer V. Lafforgue's question.}
We will now see how Proposition \ref{principale_h} and \ref{principale_o} imply that Lemmas 3.5 and 3.7 of V. Lafforgue's paper \cite{Lafforgue} hold. To see how property (RD) can be deduced we refer to the next section. But before proceeding, here are some definitions and notations. 
\begin{defi}\label{triangle_equil}Let $G$ denote $\sltk$ (for ${\bf K}={\bf R},{\bf C},{\bf H}$) or $\slto$ and $K$ denote $\sutk$ or $\suto$, and  $X_{\bf K}=G/K$. Denote by $A$ the diagonal matrices in $X_{\bf K}$. The Cartan decomposition states then that $G=K(A)$. Moreover, for any $x,y\in X_{\bf K}$, there exists $g\in G$ so that $x,y\in g(A)$. For any $t\in{\bf R}$, $x,y,z\in X_{\bf K}$, we say that $(x,y)$ is {\it of shape $(t,0)$} if there exists $g\in G$ such that
$$g(x)=I\,,\ g(y)=e^{-t/3}\left(\begin{array}{ccc}e^t&0&0\\
0&1&0\\
0&0&1\end{array}\right)$$
and that $(x,y,z)$ is an {\it equilateral triangle} of oriented size $t$ (where $t\in{\bf R}$ can be positive or negative) if there exists $g\in G$ such that
$$g(x)=I\,,\ g(y)=e^{-t/3}\left(\begin{array}{ccc}e^t&0&0\\
0&1&0\\
0&0&1\end{array}\right)\,,\ g(z)=e^{-2t/3}\left(\begin{array}{ccc}e^t&0&0\\
0&e^t&0\\
0&0&1\end{array}\right).$$

\smallskip

Let $(X,d)$ be a metric space and $\delta\geq 0$. For any finite sequence of points $x_1\dots,x_n\in X$ we say that $x_1\dots x_n$ is a {\it $\delta$-path} if
$$d(x_1,x_2)+\dots+d(x_{n-1},x_n)\leq d(x_1,x_n)+\delta$$
and that three points $x,y,z\in X$ form {\it $\delta$-retractable triple} if there exists $t\in X$ such that the paths $xty$, $ytz$ and $ztx$ are $\delta$-paths.\end{defi}
\begin{lemma}\label{13_14}
Let ${\bf K}={\bf H}$ or ${\bf O}$. For any $\delta_0>0$ there exists a $\delta>0$ such that the following is true:

1) For any $x,y,z\in X_{\bf K}$ there exists $x',y',z'$ an equilateral triangle in $X_{\bf K}$ such that the paths $xx'y'y$, $yy'z'z$, $zz'x'x $ are $\delta$-paths.

2) For any $s,t\in{\bf R}$ of same sign and $x,v,w,y\in X_{\bf K}$ such that $d(v,w)<\delta_0$, $(v,z)$ is of shape $(t,0)$ and $(w,y)$ is of shape $(s,0)$, the triangle $z,v,y$ is $\delta$-retractable.
\end{lemma}
\begin{proof} 1) Because of Proposition \ref{principale_h} and \ref{principale_o}, there is an isometric copy of $X_{\bf C}$ in $X_{\bf K}$ containing $x,y,z$. Thus, because of Lemma 3.6 in \cite{Lafforgue} there exists $\delta\geq 0$ and $x',y',z'$ an equilateral triangle in $X_{\bf C}$ such that the paths $xx'y'y$, $yy'z'z$ and $zz'x'x $ are $\delta$-paths. Since the embedding is totally geodesic, the triangle $x',y',z'$ will be equilateral in $X_{\bf K}$ as well.

\medskip

2) Without loss of generality, we can assume that $w=I$ and that $y$ is diagonal, so that $z,v\in h(A)$ for an $h\in G$ such that $d(h(I),I)\leq\delta_0$. Because of Proposition \ref{principale_h} and \ref{principale_o} there is a totally geodesic embedding $\overline{\varphi}:X_{\bf C}\to X_{\bf K}$ containing $h(I)$, $w$ and $y$. But since $\overline{\varphi}$ comes from an embedding $\varphi:\sltc\to G$ we have that:
$$h(I)=\overline{\varphi}(\overline{h})=\overline{\varphi}(h'(I))=\varphi(h')(I)$$
(for an $\overline{h}\in X_{\bf C}$ and $h'\in\sltc$) so that the whole $h(A)$ is contained in the image of $\overline{\varphi}$ in $X_{\bf K}$ and we thus can see $x,v,w$ and $y$ in $X_{\bf C}$. Using Lemma 3.7 in \cite{Lafforgue}, the triangle $z,v,y$ is $\delta$-retractable in $X_{\bf C}$, and since the embedding is totally geodesic, it will be $\delta$-retractable in $X_{\bf K}$ as well.
\end{proof}
\section{How to establish property (RD)}
We start by recalling some definitions that can be found in \cite{Lafforgue}. 
\begin{defi} Let $\delta\geq 0$. A discrete metric space $(X,d)$ satisfies {\it property $(H_\delta)$} if for there exists a polynomial $P_\delta$ such that for any $r\in{\bf R}_+$, $x,y\in X$ one has
$$\sharp\{t\in X\hbox{ such that }d(x,t)\leq r,\hbox{ and }xty\hbox{ is a }\delta-\hbox{path}\}\leq P_\delta(r).$$
We say that $(X,d)$ satisfies {\it property $(H)$} if it satisfies $(H_\delta)$ for any $\delta\geq 0$. For any $r\in{\bf R}_+$, and $x\in X$ let $B(x,r)\subset X$ denote the open ball of radius $r$ centered in $x$. A subset $Y\subset X$ is a {\it uniform net} in $X$ if there exists two constants $r_Y$ and $R_Y$ in ${\bf R}_+$ such that
\begin{eqnarray*}& B(y,r_Y)\cap Y=\{y\} & \hbox{ for any }  y\in Y\\
& B(x,R_Y)\cap Y\not=\emptyset &\hbox{ for any }  x\in X\end{eqnarray*}\end{defi}
\begin{rem}\label{H_Lie} In other words, property $(H)$ gives a polynomial bound for geodesics between two points. V. Lafforgue, in \cite{Lafforgue} proved that any uniform net in a semi-simple Lie group has property $(H)$, provided that the distance takes into account every root.
\end{rem}
\begin{defi}
Let $(X,d)$ be a metric space and $\Gamma$ be a discrete group acting by isometries on $X$. The pair $(X,\Gamma)$ satisfy {\it property $(K)$} if there exists $\delta\geq 0$, $k\in{\bf N}$ and $\Gamma$-invariant subsets ${\mathcal T}_1,\dots,{\mathcal T}_k$ of $X^3$ such that:
\begin{itemize}\item[$(K_a)$] There exists $C\in{\bf R}_+$ such that for any $(x,y,z)\in X^3$, there exists $i\in\{1,\dots,k\}$ and $(\alpha,\beta,\gamma)\in{\mathcal T}_i$ such that
$$\max\{d(\alpha,\beta),d(\beta,\gamma),d(\gamma,\alpha)\}\leq C \min\{d(x,y),d(y,z),d(z,x)\}+\delta$$
and $x\alpha\beta y$, $y\beta\gamma z$, $z\gamma\alpha x$ are $\delta$-paths.
\medskip
\item[$(K_b)$] For any $i\in\{1,\dots,k\}$ and $\alpha,\beta,\gamma,\gamma'\in X$, if $(\alpha,\beta,\gamma)\in{\mathcal T}_i$ and $(\alpha,\beta,\gamma')\in{\mathcal T}_i$ then the triangles $\alpha\gamma\gamma'$ and $\beta\gamma\gamma'$ are $\delta$-retractable.\end{itemize}
\end{defi}
\begin{thm}[V. Lafforgue \cite{Lafforgue}]\label{H+K=RD} Let $X$ be a discrete metric space, and let $\Gamma$ be a group acting freely and isometrically on $X$. If the pair $(X,\Gamma)$ satisfies $(H)$ and $(K)$, then $\Gamma$ satisfies property (RD).\end{thm}
From now on, ${\mathcal K}$ will denote a finite product of metric spaces as described in the introduction. The two following lemmas will show why any uniform net in ${\mathcal K}$ has property $(H)$.
\begin{lemma}\label{H_grossiere_invariance}Let $X$ be a uniform net in a geodesic metric space $({\mathcal H},d)$. If $X$ satisfies property $(H)$, then so will any other uniform net $Y$ in ${\mathcal H}$.\end{lemma}
\begin{proof} Take $x,y\in{\mathcal H}$ and set, for any $\delta\geq 0$,
$$\Upsilon_{(\delta,r)}(x,y)=\{t\in{\mathcal H}|xty\hbox{ is a }\delta-\hbox{path},\,d(x,t)\leq r\}.$$
Because of property $(H)$ we can cover $\Upsilon_{(\delta,r)}(x,y)$ with $P_\epsilon(r)$ balls of radius $R_X$ centered at each point of $\Upsilon_{(\delta,r)}(x,y)\cap X$, where $P_\epsilon$ denotes the polynomial associated to the constant $\epsilon=2R_X+\delta$ in the definition of property $(H)$. Now, take $z\in Y$ so that $xzy$ is a $\delta$-path. Since in particular $z\in{\mathcal H}$, we can find $t\in X$ at a distance less than $R_X$, and thus $z$ is in the ball of radius $R_X$ centered at $t$. It is now obvious that $xty$ is a $2R_X+\delta$-path. But in each ball of radius $R_X$ there is a uniformly bounded number $N$ of elements of $Y$, so that
\begin{eqnarray*}& &\sharp\{t\in Y\hbox{ such that }d(x,t)\leq r,\hbox{ and }xty\hbox{ is a }\delta-\hbox{path}\}\\
&\leq &N\sharp\{t\in X\hbox{ such that }d(x,t)\leq r,\hbox{ and }xty\hbox{ is a }(2R_X+\delta)-\hbox{path}\}\\
&\leq& N P_\epsilon(r)\end{eqnarray*}
and thus $Y$ satisfies property $(H)$, choosing $P_\delta=NP_\epsilon$.
\end{proof}
\begin{lemma}\label{surHdelta} Let $({\mathcal H}_1,d_1)\dots({\mathcal H}_n,d_n)$ be metric spaces whose uniform nets all have property $(H)$. Endow ${\mathcal H}={\mathcal H}_1\times\dots\times{\mathcal H}_n$ with the $\ell^1$ combination of the distances $d_i$, then any uniform net in ${\mathcal H}$ has property $(H)$.\end{lemma}
\begin{proof} Because of Lemma \ref{H_grossiere_invariance}, it is enough to show that one particular uniform net has property $(H)$. To do that, let $X_1\subset{\mathcal H}_1\dots X_n\subset{\mathcal H}_n$ be uniform nets and look at $X=X_1\times\dots\times X_n$, which is a uniform net in ${\mathcal H}$. Take $x,y\in X$, $\delta,r\geq 0$ and $t\in\Upsilon_{(\delta,r)}(x,y)$. We write $x,y,t$ in coordinates, that is to say $x=(x_1,\dots,x_n),y=(y_1,\dots,y_n)$ and $t=(t_1,\dots,t_n)$, where $x_i,y_i,t_i\in X_i$ for any $i=1,\dots,n$. On each factor $X_i$, $x_it_iy_i$ will be a $\delta$-path as well, and since we are considering the $\ell^1$ combination of norms, $d_i(x_i,t_i)\leq r$. On each $X_i$ there is by assumption at most $P_{\delta,i}(r)$ of those points $t_i$ and thus on $X$ we have at most
$$P_\delta(r)=\prod_{i=1}^nP_{\delta,i}(r)$$
$t$'s at distance to $x$ less than $r$ and such that $xty$ is a $\delta$-path. Obviously $P_\delta$ is again a polynomial (of degree the sum of the degrees of the $P_{\delta,i}$'s) and thus $X$ has property $(H)$.
\end{proof}
It is an easy observation that any uniform net in a Gromov hyperbolic space has property $(H)$, and thus, since any of the metric spaces forming ${\mathcal K}$ have property $(H)$, we deduce that ${\mathcal K}$ has property $(H)$ as well. Let us now see what happens to property $(K)$.
\begin{defi}\label{equil_gen} A triple $x,y,z\in{\mathcal K}$ forms an {\it equilateral triangle} if when we write in coordinates $x=(x_1,\dots,x_n)$, $y=(y_1,\dots,y_n)$ and $z=(z_1,\dots,z_n)$, the triples $x_i,y_i,z_i\in{\mathcal X}_i$ form
\begin{itemize}
\item an equilateral triangle in the sense given in Definition \ref{triangle_equil} if ${\mathcal X}_i$ is an $X_{\bf K}$.
\item an equilateral triangle in the sense given in \cite{RRS} if ${\mathcal X}_i$ is an $\widetilde{A}_2$-type building.
\item a single point (i.e. $x_i=y_i=z_i$) if ${\mathcal X}_i$ is a hyperbolic space.\end{itemize}
Since we endowed ${\mathcal K}$ with the $\ell^1$ combinations of of the norms on the ${\mathcal X}_i$'s it is obvious that in particular an equilateral triangle will satisfy $d(x,y)=d(y,z)=d(z,x)$. Concerning the orientation, remember that in Definition \ref{triangle_equil} we had two possible orientations for a triangle, positive or negative. In ${\mathcal K}$ we will have much more possible orientations since the orientation of a triangle will depend on its orientation in each non hyperbolic coordinate (in hyperbolic coordinates, the projection is just a point and thus has only one possible orientation). Suppose that among the ${\mathcal X}_i$'s forming ${\mathcal K}$, $m$ of them are not Gromov hyperbolic (for $0\leq m\leq n$), and set
$$I=\{(a_1,\dots,a_m)|a_i\in\{+,-\}\},$$
we then say that an equilateral triangle $x,y,z$ has orientation $j\in I$ if in the non hyperbolic components ${\mathcal X}_i$ the triangle $x_i,y_i,z_i$ has the orientation given by $a_i$ for $i=1,\dots m$.
\end{defi}
\begin{lemma}\label{13_14_bis}
For any $\delta_0>0$ there exists a $d>0$ such that the following is true:

1) For any $x,y,z\in{\mathcal K}$ there exists $x',y',z'$ an equilateral triangle in ${\mathcal K}$ such that the paths $xx'y'y$, $yy'z'z$ and $zz'x'x$ are $d$-paths.

2) For any two equilateral triangles $x,y,z$ and $a,b,c$ in ${\mathcal K}$ of same orientation and such that $d(x,a)$, $d(y,b)$ are both less than $\delta_0$, the triangles $x,z,c$ and $y,z,c$ are $d$-retractable.
\end{lemma}
\begin{proof}1) For any $i=1,\dots,n$ it exists an equilateral triangle $x'_i,y'_i,z'_i$ in ${\mathcal X}_i$ such that the paths $x_ix'_iy'_iy_i$, $y_iy'_iz'_iz_i$ and $z_iz'_ix'_ix_i$ are $\delta_i$-paths. Indeed, this is because of Lemma 3.6 in \cite{Lafforgue} and Lemma \ref{13_14} in case ${\mathcal X}_i$ is $X_{\bf K}$ for ${\bf K}={\bf R}, {\bf C}, {\bf H}$ or ${\bf O}$, because of Section 3 in \cite{RRS} in case ${\mathcal X}_i$ is an $\widetilde{A}_2$-type building and trivially in case ${\mathcal X}_i$ is a hyperbolic space. Now, setting $x'=(x'_1,\dots,x'_n),y'=(y'_1,\dots,y'_n)$ and $z'=(z'_i,\dots,z'_n)$ we see that by construction the triangle $x',y',z'$ is equilateral and that the path $xx'y'y$ is a $d$-paths for any $d\geq\sum_{i=1}^n\delta_i$ since
\begin{eqnarray*}d(x,x')+d(x',y')+d(y',y)&=&\sum_{i=1}^nd_i(x_i,x'_i)+d_i(x'_i,y'_i)+d_i(y'_i,y_i)\\
&\leq &\sum_{i=1}^n(d_i(x_i,y_i)+\delta_i)\leq d(x,y)+d\end{eqnarray*}
and similarly for the paths $yy'z'z$ and $zz'x'x$.

\bigbreak

2) Here we have to find a $d\geq 0$ and $u$, $v$ in ${\mathcal K}$ such that the paths $xuz,zuc$ and $cux$ as well as the paths $yvz,zvc$ and $cvy$ are $d$-paths. But for any $i=1,\dots,n$, the triangles $x_i,y_i,z_i$ and $a_i,b_i,c_i$ are equilateral triangles and thus because of Lemma 3.7 in \cite{Lafforgue} and Lemma \ref{13_14} in case ${\mathcal X}_i$ is $X_{\bf K}$ for ${\bf K}={\bf R}, {\bf C}, {\bf H}$ or ${\bf O}$, trivially in case ${\mathcal X}_i$ is an $\widetilde{A}_2$-type building or a hyperbolic space, the triangles $x_i,z_i,c_i$ and $y_i,z_i,c_i$  are $\delta_i$-retractable, that is there exists $\delta_i\geq 0$ and  points $u_i$ and $v_i$ so that the paths $x_iu_iz_i,z_iu_ic_i$ and $c_iu_ix_i$ as well as the paths $y_iv_iz_i,z_iv_ic_i$ and $c_iv_iy_i$  are $\delta_i$-paths. Now the points $u=(u_1,\dots,u_n)$ and $v=(v_1,\dots,v_n)$ are the sought points on which the triangles $x,z,c$ and $y,z,c$ retract, for $d\geq\sum_{i=1}^n\delta_i$. Indeed, let us check that for the path $xuz$:
\begin{eqnarray*}d(x,u)+d(u,z)&=&\sum_{i=1}^n(d_i(x_i,u_i)+d(u_i,z_i))\\
&\leq &\sum_{i=1}^n(d_i(x_i,z_i)+\delta_i)\leq d(x,z)+d.\end{eqnarray*}
and similarly for the paths $zuc$ and $cux$ as well as the paths $yvz,zvc$ and $cvy$.
\end{proof}
The following lemma is the analogue of a part of Theorem 3.3 in \cite{Lafforgue}.
\begin{lemma}Let $\Gamma$ be a discrete cocompact subgroup of the isometry group of ${\mathcal K}$, and $Z\subset{\mathcal K}$ be a $\Gamma$-invariant uniform net. Let $X$ be a free $\Gamma$-space and $\theta:X\to Z$ be a $\Gamma$-equivariant map. Endow $X$ with the distance 
$$\theta^*(d)(x,y)=\left\{\begin{array}{cc}0&\hbox{ if }x=y\\
1+d(\theta(x),\theta(y))&\hbox{ if }x\not=y\end{array}\right.$$ 
where $d$ is the induced distance of ${\mathcal K}$ on $Z$. Then the pair $(X,\Gamma)$ satisfies property $(K)$.\end{lemma} 
\begin{proof} Take $\delta\geq 4R_Z+d$ (for $d$ as in Lemma \ref{13_14_bis}). We first have to define the $\Gamma$-invariant subsets of $X^3$. Let us consider $I$ as explained in Definition \ref{equil_gen}, that is, $I$ is a set of indices running along the possible orientations of equilateral triangles in ${\mathcal K}$. For any $i\in I$, we define
\begin{eqnarray*}{\mathcal T}'_i=\{(\alpha,\beta,\gamma)\in Z^3&|&\hbox{ it exists }(a,b,c)\in{\mathcal K}^3\hbox{ equilateral triangle }\\
& &\hbox{ with }d(\alpha,a)\leq R_Z,d(\beta,b)\leq R_Z,d(\gamma,c)\leq R_Z\}
\end{eqnarray*}
In other words, ${\mathcal T}'_i$ is the set of triples of $Z$ which are not too far from an equilateral triangle, and since $\Gamma$ acts on ${\mathcal K}$ by isometries, the sets ${\mathcal T}'_i$ are $\Gamma$-invariant. We then set, for any $i\in I$:
$${\mathcal T}_i=\theta^{-1}({\mathcal T}'_i).$$
Since $\theta$ is $\Gamma$-equivariant, ${\mathcal T}_i$ is $\Gamma$-invariant for any $i\in I$. Let us explain why then $(X,\Gamma)$ satisfy property $(K)$. Because of the distance defined on $X$, it is enough to prove $(K_a)$ and $(K_b)$ for $Z$ and the sets ${\mathcal T}'_i$.
 
\bigskip

$\underline{(K_a)}$:

\smallskip

Take $(x,y,z)\in Z^3$, we have to show that there exists $(\alpha,\beta,\gamma)$ in some ${\mathcal T}_i$ so that the triple $(x,y,z)$ retracts on  $(\alpha,\beta,\gamma)$. But because of part a) of Lemma \ref{13_14_bis}, we know that there exists $(x',y',z')\in{\mathcal K}^3$, forming an equilateral triangle and so that the triple $x,y,z$ retracts on $x',y',z'$. Now $Z$ being a uniform net in ${\mathcal K}$, there exists $\alpha,\beta,\gamma$ three points of $Z$ with $d(\alpha,x')\leq R_Z$, $d(\beta,y')\leq R_Z$ and $d(\gamma,z')\leq R_Z$, so that the triple $(\alpha,\beta,\gamma)$ belongs to ${\mathcal T}'_i$ for some $i\in I$. We compute:
\begin{eqnarray*}& &d(x,\alpha)+d(\alpha,\beta)+d(\beta,y)\\
&\leq&d(x,x')+d(x',\alpha)+d(\alpha,x')+d(x',y')+d(y'\beta)+d(\beta,y')+d(y',y)\\
&\leq&d(x,x')+d(x',y')+d(y',y)+4R_Z\leq d(x,y)+d+4R_Z\leq\delta\end{eqnarray*} 
and similarly for the paths $y\beta\gamma z$ and $z\gamma\alpha x$.

\bigskip

$\underline{(K_b)}$:

\smallskip

Take $i\in I$ and four points in $Z$ defining two triples in ${\mathcal T}'_i$, say $(\alpha,\beta,\gamma)$ and $(\alpha,\beta,\gamma')$. We have to prove that both triangles $\alpha,\gamma,\gamma'$ and $\beta,\gamma,\gamma'$ are $\delta$-retractable. By definition of ${\mathcal T}'_i$ we can find two equilateral triangles $a,b,c$ and $x,y,z$ such that 
$$d(\alpha,a)\leq R_Z,d(\beta,b)\leq R_Z,d(\gamma',c)\leq R_Z$$ 
and 
$$d(\alpha,x)\leq R_Z,d(\beta,y)\leq R_Z,d(\gamma,z)\leq R_Z$$ 
so that obviously $d(a,x)\leq 2R_Z$ and $d(b,y)\leq 2R_Z$ and thus applying part b) of Lemma \ref{13_14_bis} we have the existence of a $d\geq 0$, and two points $u$ and $v$ in ${\mathcal K}$ so that the paths $xuz,zuc$ and $cux$ as well as the paths $yvz,zvc$ and $cvy$ are $d$-paths. Again, $Z$ being a uniform net in ${\mathcal K}$, we can find $u'$ and $v'$ in $Z$ at respective distances less than $R_Z$ to $u$ and $v$. We claim that the paths $\alpha u'\gamma$, $\gamma u'\gamma'$ and $\gamma'u'\alpha$ as well as the paths $\beta v'\gamma$, $\gamma v'\gamma'$ and $\gamma'v'\beta$ are $\delta$-paths:
\begin{eqnarray*}& &d(\alpha,u')+d(u',\gamma)\\
&\leq&d(\alpha,x)+d(x,u)+d(u,u')+d(u',u)+d(u,z)+d(z,\gamma)\\
&\leq&4R_Z+d(x,u)+d(u,z)\leq d(x,z)+d+4R_Z\\
&\leq& d(x,\alpha)+d(\alpha,\gamma)+d(\gamma,z)+d+4R_Z\\
&\leq& d(\alpha,\gamma)+(d+6R_Z)\leq d(\alpha,\gamma)+\delta\end{eqnarray*}
and similarly for the other paths.
\end{proof}
Now, if under the assumptions of this lemma we furthermore assume that $\sharp\theta^{-1}(z)\leq N$ (i. e. $\theta$ has uniformly bounded fibers) we use Remark \ref{H_Lie} and Lemma \ref{H_grossiere_invariance} to deduce that $X$ has property $(H)$. We now can apply Theorem \ref{H+K=RD}, as follows: $Z$ is a $\Gamma$-invariant uniform net in ${\mathcal K}$ and let $Z=\coprod_{j\in J}\Gamma x_j$ its partition in $\Gamma$-orbits. Then with $X=\coprod_{i\in J}\Gamma$ and $\theta$ the obvious orbit map and with $\{{\mathcal T}_i\}_{i\in I}$ as defined in the previous lemma we get:
\begin{thm} Let $\Gamma$ be a discrete group acting by isometries on ${\mathcal K}$ and with uniformly bounded stabilizers on some $\Gamma$-invariant uniform net. Then $\Gamma$ has property (RD).\end{thm}
\hfill\qedsymbol

If $\Gamma$ is a cocompact lattice of isometries on ${\mathcal K}$, it is enough to take $Z=\Gamma x_0$ and $X=\Gamma$, so $\Gamma$ has property (RD).
\begin{rem} If one is only interested in cocompact lattices in the isometry group of a product of hyperbolic spaces the proof of property (RD) then becomes much simpler. Indeed, denote by ${\mathcal H}$ a finite product of hyperbolic spaces, $\Gamma$ a cocompact lattice in its isometry group and $X=\Gamma x_0$ for some $x_0\in{\mathcal H}$. We showed that $X$ has property $(H)$, so it remains to show that any triple of points is $\delta$-retractable, but this is a special case of Lemma \ref{13_14_bis} part 1). We can conclude using Proposition 2.3 of \cite{Lafforgue}.

Notice that for this work we have established property (RD) for $\Gamma$ a cocompact lattice of isometries on ${\mathcal K}$, endowed with the $\ell^1$ combination of the distances, and this implies Theorem \ref{principal} since ${\rm Iso}({\mathcal X}_1)\times\dots\times{\rm Iso}({\mathcal X_n})\subset{\rm Iso}({\mathcal K},\ell^1)$. The arguments used for the proof fail for other combinations of the distances.\end{rem}
%

\end{document}